\documentclass[10pt]{article}

\usepackage{amsmath,amscd,amsfonts,amsthm}

\DeclareMathOperator{\rad}{rad}

\renewcommand{\aa}{\mathbf{s}}
\newcommand{\ep}{\epsilon}
\newcommand{\beq}{\begin{displaymath}}
\newcommand{\eeq}{\end{displaymath}}

\theoremstyle{plain}
\newtheorem{thm}{Theorem}

\newtheorem{lem}[thm]{Lemma}

\theoremstyle{definition}

\newtheorem{conj}[thm]{Conjecture}

\theoremstyle{remark}

\title{Congruence ABC implies ABC}
\author{Jordan S. Ellenberg \\ Princeton University \\
\texttt{ellenber@math.princeton.edu}}
\date{25 March 1999}

\begin{document}
\maketitle

\begin{abstract}
The ABC conjecture of Masser and Oesterl\'{e} states that if $(a,b,c)$
are coprime integers with $a + b + c = 0$, then $\sup(|a|,|b|,|c|) <
c_\ep(\rad(abc))^{1+\ep}$ for any $\ep > 0$.  In \cite{oest:fermat},
Oesterl\'{e} observes that if the ABC conjecture holds for all
$(a,b,c)$ with $16 | abc$, then the full ABC conjecture holds.  We
extend that result to show that, for every integer $N$, the
``congruence ABC conjecture'' that ABC holds for all $(a,b,c)$ with
$N|abc$ implies the full ABC conjecture.
\end{abstract}

\section{Introduction}

The ABC conjecture was introduced by Masser and Oesterl\'{e} in 1985,
and has since been shown to be related to many other conjectures,
especially conjectures regarding the arithmetic of elliptic
curves~\cite{mazu:abc}. 

For our purposes, an {\em ABC-solution} $\aa$ is a triple $(a,b,c)$ of distinct
relatively prime integers satisfying $a + b + c = 0$, and such that
$a$ and $b$ are negative.  (The requirement that the integers be
distinct is included only to simplify the exposition below.)  If $n>0$
is an integer, the radical $\rad(n)$ is defined to be the product of
all primes dividing $n$.

For any $\ep>0$, we define a function on ABC-solutions by
\beq
f(\aa,\ep) = \log(c) - (1+\ep)\log \rad(abc).
\eeq

Then the ABC conjecture can be phrased as follows:

\begin{conj}[ABC conjecture]
For each $\ep > 0$, there exists a constant $C_\ep$ such that 
\beq
f(\aa,\ep) < C_\ep
\eeq
for all $\aa$.
\end{conj}

In \cite{oest:fermat}, Oesterl\'{e} showed that the ABC conjecture is
equivalent to a conjecture of Szpiro on elliptic curves (\cite[Conj.\
4]{oest:fermat}) In the proof, he observes that if the ABC conjecture
is known to hold for all $(a,b,c)$ with $16 | abc$, then the full ABC
conjecture can be shown to hold.  This suggests considering a family
of weaker conjectures indexed by integers $N$, as follows:

\begin{conj}[Congruence ABC conjecture for $N$]
For each $\ep > 0$, there exists a constant $C_\ep$ such that
\beq
f(\aa,\ep) < C_\ep
\eeq
for all $\aa$ such that $N | abc$.
\end{conj}

It has long been known to experts that the congruence ABC conjecture
for any $N$ is equivalent to the full ABC conjecture.  However, a
proof has never to our knowledge appeared in the literature, and we
take the opportunity to provide one in this note.

\section{Congruence ABC implies ABC}

\begin{thm}
The congruence ABC conjecture for $N$ implies the ABC conjecture.
\label{th:congabc}
\end{thm}

\begin{proof}
For each positive even integer $n$, we define an operation $\Theta_n$ on
ABC-solutions as follows.  Let $\aa = (a,b,c)$ be an ABC-solution.
Then
\beq
\Theta_n(\aa) = (-2^{-m}(a-b)^n, -2^{-m}[c^n - (a-b)^n], 2^{-m}c^n)
\eeq
where $m = n$ if $c$ is even, and $m = 0$ otherwise.  Then
$\Theta_n(\aa)$ is again an ABC-solution.

\begin{lem}
There exist constants $c_{n,\ep}>0$ and $c'_{n,\ep}$ such that 
\beq
f(\Theta_n(\aa),\ep/[n + (n-1)\ep]) \geq c_{n,\ep}f(\aa,\ep) + c'_{n,\ep}.
\eeq
\end{lem}

\begin{proof}
Let $A = -2^{-m}(a-b)^n, B = 2^{-m}(-c^n + (a-b)^n), C = 2^{-m}c^n$.  Then
 
\beq
\log \rad(ABC) \leq
\log|a-b| + \log \rad(abc) + \log \rad(B/ab).
\eeq

Now

\begin{eqnarray*}
B/ab & = & \frac{(a+b)^n - (a-b)^n}{ab} \\
& = & 4[(a+b)^{n-2} + (a+b)^{n-4}(a-b)^2 + \ldots + (a-b)^{n-2}] \\
& \leq & 2n(a+b)^{n-2}.
\end{eqnarray*}

So

\begin{eqnarray*}
\log \rad(ABC) & \leq & \log|a-b| + \log \rad(abc) + (n-2)\log c +
\log 2n \\
& \leq & (n-1) \log c + \log \rad (abc) + \log 2n \\
& = & (n-1) \log c + (1 + \ep)^{-1}(\log c - f(\aa,\ep)) + \log 2n \\
&  = & n \log c - \ep(1 + \ep)^{-1} \log c  - (1+\ep)^{-1}
f(\aa,\ep) + \log 2n \\
& = & (1 - \ep[n(1 + \ep)]^{-1})(\log C + m \log 2) \\
& &    - (1+\ep)^{-1}f(\aa,\ep) + \log 2n. \\
\end{eqnarray*}

It follows that
\beq
\log C \geq \frac{n(1+\ep)}{(n + n\ep - \ep)} (\log \rad(ABC) + (1 + \ep)^{-1}f(\aa,\ep) -
\log 2n) - n \log 2.
\eeq

Thus

\beq
f(\Theta_n(\aa),\ep/[n+n\ep-\ep]) \geq c_{n,\ep}f(\aa,\ep) + c'_{n,\ep}.
\eeq

where
\beq
c_{n,\ep} = n/(n + n\ep - \ep)
\eeq
and
\beq
c'_{n,\ep} = -n(1+\ep)(\log 2n)/(n + n\ep - \ep) - n \log 2.
\eeq
\end{proof}

We now proceed with the proof of Theorem~\ref{th:congabc}.  Assume
that the congruence ABC conjecture for $N$ is true.  Then there exists
$C_\ep$ such that 
\beq
f(\aa,\ep) < C_\ep
\eeq
for all $\aa$ with $N | abc$.

Let $n = \phi(N)$.  If $N=2$, Theorem~\ref{th:congabc} is trivial.  We
may therefore assume that $n$ is even.

\begin{lem} 
Let $(A,B,C) = \Theta_n(\aa)$.  Then $N | ABC.$
\end{lem}

\begin{proof}
Suppose $p$ is an odd prime dividing $N$, and let $p^\nu$ be the
largest power of $p$ dividing $N$.  Then $(p-1)p^{\nu-1} | n$.  In
particular, $\nu < n$.  If $p$ divides $c$ or $a-b$, then $p^{\nu} |
p^n |ABC$.  If $p$ divides neither $c$ nor $a-b$, then $-A
=2^{-m}(a-b)^n$ and $C = 2^{-m}c^n$ are{\ c}ongruent mod $p^\nu$;
therefore, $p^\nu | B$.

Now let $2^\nu$ be the largest power of $2$ dividing $N$.  If $c$ is
even, then so is $a-b$, and exactly one of $c$ and $a-b$ is a multiple
of $4$.  Thus, one of $(a-b)^n$ and $c^n$ is a multiple of $4^n$.
Since, in this case, $A = 2^{-n}(a-b)^{n}$ and $C = 2^{-n}c^n$, one of
$A$ and $C$ is a multiple of $2^n$, whence also of $2^\nu$.  If, on
the other hand, $c$ is odd, then so is $a-b$.  Then $(a-b)^n$ and
$c^n$ are both congruent to $1$ mod $2^nu$, so $2^\nu | B$.
\end{proof}

For any $\ep > 0$, it now follows that
\beq
f(\aa,\ep) \leq c_{n,\ep}^{-1}(f(\Theta_n(\aa),\ep/[n + n\ep - \ep]) - c'_{n,\ep})
\leq c_{n,\ep}^{-1}(C_{\ep/[n + n\ep - \ep]} - c'_{n,\ep})
\eeq
for any ABC-solution $\aa$.  Since the right hand side depends only on
$\epsilon$ and $N$, this proves the full ABC conjecture.
\end{proof}

\section*{Acknowledgments}
It is my pleasure to thank Joseph Oesterl\'{e} and the referee for
helpful comments on this note.

\end{document}